\documentclass[12pt]{amsart}

\usepackage{amsfonts,amsmath,amsthm}
\usepackage{latexsym}

\newtheorem{theorem}{Theorem}[section]

\newtheorem{lemma}[theorem]{Lemma}

\newtheorem{definition}{Definition}[section]
\theoremstyle{definition}

\newcommand \tr {\mathrm{Tr}\:}

\newcommand{\C}{\ensuremath{\mathbb{C}}}
\newcommand{\R}{\ensuremath{\mathbb{R}}}
\newcommand{\K}{\ensuremath{\mathbb{K}}}
\newcommand{\adj}{^\ast}

\def \l {\lambda}

\def \< {\left\langle}
\def \> {\right\rangle}

\begin{document}

\title{A geometric estimate on the norm of product of functionals}

\maketitle
\centerline{\Large{M\'at\'e Matolcsi\footnote{Author supported by 
Hungarian research funds OTKA-T047276, 
OTKA-F049457, OTKA-T049301.}}}
\medskip
\centerline{Alfr\'ed R\'enyi Institute of Mathematics,}
\centerline{Hungarian Academy of Sciences,}
\centerline{POB 127, H-1364, Budapest, Hungary.}
\centerline{Tel: (+361) 483-8302, Fax: (+361) 483-8333,}
\centerline{e-mail: matomate@renyi.hu}

\begin{abstract}
The open problem of determining the exact value of the $n$-th 
linear polarization constant $c_n$ of $\R^n$ has received 
considerable attention over the past few years.  
This paper makes a contribution to the subject by providing a new lower bound on the value of 
 $\sup_{\|{\bf{y}}\|=1}|\< {\bf{x}}_1,{\bf{y}} \> \cdots \< {\bf{x}}_n,{\bf{y}} \> |$, where ${\bf{x}}_1,
\dots ,{\bf{x}}_n$ are unit vectors in $\R^n$. The new estimate is given in terms of the eigenvalues of the Gram matrix $[\< {\bf{x}}_i,{\bf{x}}_j\>  ]$ and 
improves upon earlier estimates of this kind. However, the intriguing
conjecture $c_n=n^{n/2}$ remains open.
\end{abstract}



\section{Introduction}
The present work contributes to study of the $n$-th linear
polarization constant $c_n(\R^n)$ of the $n$-dimensional real Euclidean space. We
begin with introducing some (more general) standard terminology
and giving a short account of some related results.

Let $X$ denote a Banach space over the real or complex field $\K$. A function
$P: \ X\to \K$ is a continuous {\it n-homogeneous polynomial} if there exists
a symmetric, continuous $n$-linear form $L: \ X^n\to\K$ such that
$P({\bf{x}})=L({\bf{x}},\dots ,{\bf{x}})$ for all ${\bf{x}}\in X$. We define
$$\|P\|:=\sup \{ |P({\bf{x}})|: \ {\bf{x}}\in B\}$$
where $B$ denotes the unit ball of $X$. Considerable attention has been
devoted to polynomials of the form $P({\bf{x}})=f_1({\bf{x}})f_2({\bf{x}})\dots f_n({\bf{x}})$, where
$f_1, f_2,\dots ,f_n$ are bounded linear functionals on $X$. For any{\it complex}
Banach space $X$ Ben\'\i tez, Sarantopoulos and Tonge \cite{ben} have obtained
$$\|f_1\| \ \|f_2\|\cdots \|f_n\|\le n^n\|f_1 f_2 \dots f_n\|,$$
and they also showed that, in general, the constant $n^n$ is best possible. For {\it real} Banach spaces, Ball's
solution \cite{ball2} of the famous plank problem of Tarski gives the same result. For specific spaces, however, the general constant $n^n$ can
be lowered. This fact motivated the following  

\begin{definition} (Ben\'\i tez, Sarantopoulos, Tonge \cite{ben})
The n-th linear polarization constant of a Banach space $X$ is defined by
\begin{eqnarray*}\label{polconstdef}
c_n(X)&:=& \inf \{M : \|f_{1}\| \cdots \|f_{n}\| \leq M
\|f_{1}\cdots f_{n}\|\; (\forall f_{1}, \ldots ,f_{n} \in
X^{\ast}) \}
 \\
 & = & 1/ \inf_{f_{1}, \ldots ,f_{n} \in S_{X^{\ast}}}
\sup_{\|{\bf{x}}\|=1}|f_{1}({\bf{x}}) \cdots f_{n}({\bf{x}})|.
\end{eqnarray*}
The linear polarization constant of $X$ is defined by
\begin{equation}\label{limitconstant}
c(X):=\lim_{n \rightarrow \infty} c_{n} (X)^{ \frac{1}{n}}\;\;.
\end{equation}
\end{definition}

Let us recall that the above definition of $c(X)$ is justified since
R\'ev\'esz and Sarantopoulos \cite{sar} showed that the limit
\eqref{limitconstant} does exist. Moreover, they also showed (both
in the real and complex cases) that $c(X)=\infty$ if and only if
$\dim X=\infty$.

Note that it is easy to see that for any Banach space $X$ we have
\begin{equation}\label{finite}
c_{n}(X)=\sup \left\{c_{n} (Y): Y {\mbox {is a subspace
of}}\, X, \, \dim Y=n \right\}\,.
\end{equation}
In particular, for a real or complex Hilbert space $H$
of dimension at least $n$, we
always have $c_n(H)=c_n(\K^n)$.

Ben\'\i tez, Sarantopoulos and Tonge \cite{ben} proved that for
isomorphic Banach spaces $X$ and $Y$ we have $c_n(X)\le
d^n(X,Y)c_n(Y)$, where $d(X,Y)$ denotes the Banach-Mazur distance
of $X$ and $Y$. Note, that for any $n$-dimensional space $X$ a
result of John \cite{john} states that $d(X,\K^n)\le \sqrt{n}$
(where $\K^n$ denotes the $n$-dimensional Hilbert space).  The
combination of these results mean that the determination of
$c_n(\K^n)$ gives information on the linear polarization constants
of other spaces, too.

In this paper we are going to focus our attention to Hilbert spaces.
Pappas and R\'ev\'esz \cite{papp} showed that $c(\K^n)= e^{-L(n,\K
)}$, where
$$L(n,\K ):= \int_{S}\textrm{log}|\< {\bf{x}},{\bf{e}}\>  |d\sigma ({\bf{x}});$$
here $S$ and $\sigma$  denote the unit sphere and the normalized
surface measure, respectively, and ${\bf{e}}\in S$ is an arbitrary unit
vector. This result gives information on the asymptotic behaviour
of $c_m(\K^n)$ as $m\to \infty$.
However, the exact values of $c_m(\K^n)$ seem, in general, hopeless  to determine. 

A remarkable result of
Arias-de-Reyna \cite{rey} states that $c_n(\C^n)=n^{n/2}$. Ball's
recent solution \cite{ball} of the complex plank problem also
implies the same result.

Compared to the complex case, the value of $c_n(\R^n)$ seems harder to find. The
determination of $c_n(\R^n)$, by the definition and the Riesz
representation theorem, boils down to determining
$$
I:=
\inf_{{\bf{x}}_{1}, \dots ,{\bf{x}}_{n} \in S} \sup_{\|{\bf{y}}\|=1}|\< {\bf{x}}_{1},{\bf{y}}\>
\cdots \< {\bf{x}}_{n},{\bf{y}}\>  |
$$
The estimate $I\le n^{-\frac{n}{2}}$ follows by considering an
orthonormal system.

The complex result of Arias-de-Reyna can be used to derive the following
estimates (see \cite{sar}, where the argument is based on an interesting complexification result of \cite{munoz}):
\begin{equation}\label{legjobb}
n^{\frac{n}{2}}\le c_n(\R^n)\le 2^{\frac{n}{2}-1}n^{\frac{n}{2}}.
\end{equation}
A natural, intriguing conjecture, see \cite{ben}, \cite{sar} is
the following.

\medskip
\noindent {\bf Conjecture.} $c_n(\R^n)=n^{n/2}$.
\medskip

Marcus (communicated in \cite{marc}, and elaborated later in
\cite{sar}) gives the following estimate: If ${\bf{x}}_1, {\bf{x}}_2, \dots ,
{\bf{x}}_n$ are unit vectors in $\R^n$ then there exists a unit vector
${\bf{y}}$ such that
\begin{equation}\label{marcus}
|\< {\bf{x}}_1,{\bf{y}}\> \cdots\< {\bf{x}}_n,{\bf{y}}\> |\ge (\l_1/n)^{n/2},
\end{equation}
where $\l_1$ denotes the smallest eigenvalue of the Gram matrix
$XX\adj =[\< {\bf{x}}_i,{\bf{x}}_j\>  ]$. Marcus also expressed the opinion
that lower bounds on $\sup_{\|{\bf{y}}\|=1}|\< {\bf{x}}_1,{\bf{y}}\>  \cdots \<
{\bf{x}}_n,{\bf{y}}\>  |$ should involve
the eigenvalues $\l_1,\dots,\l_n$ of the Gram matrix  $XX\adj
=[\< {\bf{x}}_i,{\bf{x}}_j\>  ]$, i.e. we should look for estimates of the form
$\sup_{\|{\bf{y}}\|=1}|\< {\bf{x}}_1,{\bf{y}}\>  \cdots \< {\bf{x}}_n,{\bf{y}}\>  |\ge f(\l_1, \dots
, \l_n) n^{-n/2}$. Note that $\sum_j \l_j=\tr XX\adj=n$.
Therefore the above Conjecture can be formulated as
$$
\sup_{\|{\bf{y}}\|=1}|\< {\bf{x}}_1,{\bf{y}}\>  \cdots \< {\bf{x}}_n,{\bf{y}}\>  | \ge 1\cdot
n^{-n/2}= \left(\frac{\l_1 +\dots +\l_n}{n}\right)^{n/2}n^{-n/2}.
$$

In \cite{matpolar} the author proved that Marcus' estimate \eqref{marcus} can be improved to 
$$
\sup_{\|{\bf{y}}\|=1}|\< {\bf{x}}_1,{\bf{y}}\>  \cdots \< {\bf{x}}_n,{\bf{y}}\>  | \ge  \left( \frac{n}{\l_1^{-1}+ \cdots +\l_n^{-1}} \right)^{n/2}n^{-n/2}.
$$
In the next section we will improve this result by replacing the harmonic mean of the numbers $\l_1, \dots , \l_n$ by the geometric mean. Also, in the course of the proof we use two 'geometrical' lemmas which may be of independent interest. The original Conjecture (invoving the arithmethic mean of the numbers $\l_1, \dots , \l_n$), however, still remains open.

\section{A geometric lower bound}

For the sake of simplicity we introduce the following notations:

Let $b_n$ denote the volume of the $n$-dimensional closed unit ball $B^{n}$ (we will not need the explicit value of $b_n$). Also, let $H_\alpha :=\{{\bf{z}}=(z_1, \dots , z_n)\in \R^n: \ |\prod_{j=1}^nz_j|\ge \alpha \cdot n^{-n/2}\}$. 

In order to prove our main result, Theorem \ref{theo1},  we will need the following two geometrical lemmas:

\begin{lemma}\label{lem1}
Let $E$ be an $n$-dimensional ellipsoid symmetric with respect to the origin (i.e. the image of the $n$-dimensional unit ball under a linear transformation of full rank) of volume $Vb_n$. Assume that the $n-1$-dimensional 'horizontal slice' $E_0:=\{{\bf{z}}=(z_1, z_2, \dots z_n)\in E : \ z_n=0\}$ has $n-1$-dimensional volume $Sb_{n-1}$. Then the horizontal slice at height $h$, $E_h:=\{{\bf{z}}=(z_1, z_2, \dots z_n)\in E : \ z_n=h\}$ has $n-1$-dimensional volume 
$$f(V,S,h)=\left \{ \begin{array}{ll} 
(1-(\frac{S}{V}h)^2)^{\frac{n-1}{2}}b_{n-1} & \mathrm{if} \  |h|\le V/S\\
0 & \mathrm{if}  \ |h|> V/S
\end{array}\right.
$$
\end{lemma}
\begin{proof}
The essence of the lemma is that the function $f$ depends only on $V,S$ and $h$ and not on the actual 'shape' of the ellipsoid. 

The statement of the lemma is clear if $E$ is a 'circular ellipsoid' whose axes are the same as the coordinate axes, i.e. $E$ is the image of the unit ball $B^{n}$ under the diagonal transformation 
$$T:=\left ( 
\begin{array} {cccccc}
S^{1/(n-1)}& \ & \ & \ & \ & 0 \\
 \ & . & \ & \ & \ & \ \\
 \ & \ & . & \ & \ & \ \\
 \ & \ & \ & . & \ & \ \\
 \ & \ & \ & \ & S^{1/(n-1)} & \ \\
 0 & \ & \ & \ & \ & \frac{V}{S} 
\end{array}
\right)   
$$

In the general case, let $E=A[B^{n}]$ be the image of the unit ball $B^{n}$ under some  transformation $A$, and assume that it posesses the    prescribed parameters $V, \ S$, and let the height $h$ also be given. 
The natural idea of the proof is that we transform the ellipsoid $E$ to a circular ellipsoid whose axes are the coordinate axes and whose parameters are the same.

Let ${\bf{r}}:=(r_1, \dots ,r_n)\in E$ denote the point of $E$ whose last coordinate $r_n$ is maximal among the points of $E$, and let ${\bf{q}}\in B^{n}$ be its inverse image, i.e. ${\bf{q}}=A^{-1}{\bf{r}}$. Let $L_0:=A^{-1}[E_0]$. 

Note that the $n-1$-dimensional $L_0$ is orthogonal to the vector ${\bf{q}}$, therefore there exists a unitary transformation $U$ which takes the horizontal slice $B_0^{n}$ of $B^{n}$ to $L_0$ and the vertical unit vector ${\bf{e}}_n$ to ${\bf{q}}$ (note that if $n-1>1$ then $U$ is not uniquely determined). Then we have $AU[B^{n}]=E$, $AU[B_0^{n}]=E_0$ and $AU({\bf{e}}_n)={\bf{r}}$. 

The transformation $AU$ maps the horizontal $n-1$-dimensional hyperplane $P_0:=\{{\bf{x}}=(x_1, \dots ,x_n)\in \R^n : \ x_n=0\}$ onto itself. Denote the restriction of $AU$ to $P_0$ by $C_0$.  
Now, take the $n-1$-dimensional transformation $\tilde{C}_1:=S^{1/(n-1)}C_0^{-1}$ of the horizontal hyperplane $P_0$. This preserves $n-1$-dimensional volume (i.e. it has determinant $\pm 1$) and takes the horizontal slice $E_0$ to  $S^{1/(n-1)}B_0^{n}$.

Consider now the $n$-dimensional transformation 
$$ C_1:=\left ( 
\begin{array} {cc}
\tilde{C}_1&{\bf{0}}\\
{\bf{0}}^T&1
\end{array}
\right) .   
$$

It is clear that applying this transformation to $E$ the image ellipsoid $C_1[E]$ will posses the same parameters as $E$, i.e. the same volume $V$, the same $n-1$-dimensional volume $S$ of its horizontal slice $(C_1[E])_0=S^{1/n-1}B_0^{n}$, and the image $C_1[E_h]$  will still be at height $h$ and have the same $n-1$-dimensional volume as $E_h$.

Let ${\bf{s}}:=(s_1, \dots ,s_n)=C_1{\bf{r}}$. 

Next we consider the transformation 
$$C_2:=\left ( 
\begin{array} {cccccc}
 1 & \ & \ & \ & 0 & -s_1/s_n \\
 \ & . & \ & \ & \ & . \\
 \ & \ & . & \ & \ & . \\
 \ & \ & \ & . & \ & . \\
 \ & \ & \ & \ & 1 & -s_{n-1}/s_n \\
 0 & \ & \ & \ & 0 & 1
\end{array}
\right) 
$$

Once again it is clear that the image $C_2C_1[E]$ has the same parameters $V, \ S$ as $E$, and $(C_2C_1[E])_h=C_2C_1[E_h]$ with equal $n-1$-dimensional volume.  
To finish the proof it is enough to observe that $C_2C_1[E]=C_2C_1AU[B^{n}]=T[B^{n}]$ with the diagonal transformation $T$ above.  
\end{proof}

The next lemma establishes the connection between ellipsoids and products of functionals. 

\begin{lemma}\label{lem2}
Assume $E$ is an $n$-dimensional ellipsoid of volume $Vb_n$ (not necessarily centered at the origin). Then $E\cap H_V\ne \emptyset$. 
\end{lemma}
\begin{proof}
The proof proceeds by induction with respect to $n$. 

For $n=1$ the statement is clear.

For an arbitrary $n$ 
let ${\bf{c}}:=(c_1, \dots c_n)$ denote the centre of $E$, and assume, without loss of generality, that $c_n\ge 0$. Let $Sb_{n-1}$ be the $n-1$ dimensional volume of the horizontal slice $E\cap P_{c_n}$, where $P_{c_n}:=\{{\bf{x}}=(x_1, \dots ,x_n)\in \R^n : \ x_n=c_n\}$.  
Now, let $h:=\frac{V}{S\sqrt{n}}$ and consider the horizontal hyperplane $P:=P_{c_n+h}$. $P$ is an $n-1$-dimensional space, and $P\cap H_V=\{{\bf{z}}=(z_1, \dots z_{n-1}): \ |\prod_{j=1}^{n-1}z_j|\ge \frac{V}{c_n+h}n^{-n/2}\}\supset \{{\bf{z}}=(z_1, \dots z_{n-1}): \ |\prod_{j=1}^{n-1}z_j|\ge \frac{V}{h}n^{-n/2}\}$. Furthermore, the $n-1$-dimensional volume of $P\cap E$ is $(1-\frac{1}{n})^\frac{n-1}{2}Sb_{n-1}$ in view of Lemma \ref{lem1} and the choice of $h$. 

Finally, observe that $\frac{V}{h}n^{-n/2}= (n-1)^\frac{-n+1}{2}\left( \frac{V}{h}(1-\frac{1}{n})^\frac{n-1}{2}\frac{1}{\sqrt{n}}\right )$ and 
$(1-\frac{1}{n})^\frac{n-1}{2}S= \frac{V}{h}(1-\frac{1}{n})^\frac{n-1}{2}\frac{1}{\sqrt{n}}$, therefore the inductive hypothesis applies.
\end{proof}

We are now in position to prove our new estimate on the norm of product of functionals. 

\begin{theorem}\label{theo1}
Let unit vectors ${\bf{x}}_1,\dots , {\bf{x}}_n$ be given in $\R^n$, and 
let $\l_1, \dots ,\l_n$ denote the eigenvalues of the Gram matrix $XX\adj =[\< {\bf{x}}_i,{\bf{x}}_j\> ]$ (the matrix $X$ is formed by the given vectors as rows). Then
\begin{equation}\label{geo}
\sup_{\|{\bf{y}}\|=1}|\< {\bf{x}}_1,{\bf{y}}\>  \cdots \< {\bf{x}}_n,{\bf{y}}\>  | \ge
\left (\prod_{j=1}^n\l_j\right )^{1/2}\cdot
n^{-n/2}
\end{equation}
\end{theorem}

\begin{proof}
We may assume that the vectors ${\bf{x}}_1, {\bf{x}}_2, \dots {\bf{x}}_n$ are linearly independent, otherwise the right hand side of the inequality is 0, and the estimate is meaningless. (We remark that other considerations , such as the ones in \cite{matpolar},  also show that if we find a way to prove a good estimate in the case of linearly dependent vectors then we may get close to proving the original Conjecture. However, at present, there seems to be no better estimate than \eqref{legjobb} for the linearly dependent case.) 

The image $E$ of the unit ball $B^{n}$ under the transformation $X$ is an $n$-dimensional ellipsoid of volume $V=\left (\prod_{j=1}^n\l_j\right )^{1/2}b_n$, therefore Lemma \ref{lem2} gives the required result. 

\end{proof}

Finally, let us make the following remarks.

An advantage of the proof applied above is that it is constructive in the sense that following the constructions of Lemma \ref{lem2} we can actually find a vector ${\bf{y}}$ which satisfies \eqref{geo}.
It is clear, however, that the estimate \eqref{geo} does not settle the original 
Conjecture.

\end{document}